\documentclass[10pt,draft]{amsart}
\usepackage{amsmath,amssymb,amsthm}
\begin{document}
\newtheorem{lem}{Lemma}[section]
\newtheorem{prop}{Proposition}[section]
\newtheorem{cor}{Corollary}[section]
\numberwithin{equation}{section}
\newtheorem{thm}{Theorem}[section]
\theoremstyle{remark}
\newtheorem{example}{Example}[section]
\newtheorem*{ack}{Acknowledgment}
\theoremstyle{definition}
\newtheorem{definition}{Definition}[section]
\theoremstyle{remark}
\newtheorem*{notation}{Notation}
\theoremstyle{remark}
\newtheorem{remark}{Remark}[section]
\newenvironment{Abstract}
{\begin{center}\textbf{\footnotesize{Abstract}}%
\end{center} \begin{quote}\begin{footnotesize}}
{\end{footnotesize}\end{quote}\bigskip}
\newenvironment{nome}
{\begin{center}\textbf{{}}%
\end{center} \begin{quote}\end{quote}\bigskip}

\newcommand{\triple}[1]{{|\!|\!|#1|\!|\!|}}
\newcommand{\xx}{\langle x\rangle}
\newcommand{\ep}{\varepsilon}
\newcommand{\al}{\alpha}
\newcommand{\be}{\beta}
\newcommand{\de}{\partial}
\newcommand{\la}{\lambda}
\newcommand{\La}{\Lambda}
\newcommand{\ga}{\gamma}
\newcommand{\del}{\delta}
\newcommand{\Del}{\Delta}
\newcommand{\sig}{\sigma}
\newcommand{\ome}{\omega}
\newcommand{\Ome}{\Omega}
\newcommand{\C}{{\mathbb C}}
\newcommand{\N}{{\mathbb N}}
\newcommand{\Z}{{\mathbb Z}}
\newcommand{\R}{{\mathbb R}}
\newcommand{\Rn}{{\mathbb R}^{n}}
\newcommand{\Rnu}{{\mathbb R}^{n+1}_{+}}
\newcommand{\Cn}{{\mathbb C}^{n}}
\newcommand{\spt}{\,\mathrm{supp}\,}
\newcommand{\Lin}{\mathcal{L}}
\newcommand{\SSS}{\mathcal{S}}
\newcommand{\F}{\mathcal{F}}
\newcommand{\xxi}{\langle\xi\rangle}
\newcommand{\eei}{\langle\eta\rangle}
\newcommand{\xei}{\langle\xi-\eta\rangle}
\newcommand{\yy}{\langle y\rangle}
\newcommand{\dint}{\int\!\!\int}
\newcommand{\hatp}{\widehat\psi}
\renewcommand{\Re}{\;\mathrm{Re}\;}
\renewcommand{\Im}{\;\mathrm{Im}\;}

\title
{{On the Yamabe equation with rough potentials}}
\author{}

\author{Francesca Prinari and Nicola Visciglia}

\address{Francesca Prinari\\
Dipartimento di Matematica Universit\`a di Lecce\\
Via Provinciale Lecce--Arnesano, 73100 Lecce, Italy}

\email{francesca.prinari@unile.it}

\address{Nicola Visciglia\\
Dipartimento di Matematica Universit\`a di Pisa\\
Largo B. Pontecorvo 5, 56100 Pisa, Italy}

\email{viscigli@dm.unipi.it}


\maketitle
\date{}

\begin{abstract}

We study the existence of non--trivial solutions
to the Yamabe equation:
$$-\Delta u+ a(x)= \mu u|u|^\frac4{n-2} \hbox{ } \mu >0, x\in \Omega \subset {\mathbf R}^n
\hbox{ with }  n\geq 4,$$
$$ u(x)=0 \hbox{ on } \partial \Omega$$
under weak regularity assumptions on the potential $a(x)$.

\noindent More precisely in dimension $n\geq 5$ we assume that:
\begin{enumerate}
\item $a(x)$ belongs to the Lorentz space $L^{\frac n2, d}(\Omega)$ 
for some $1\leq d
<\infty$,
\item
 $a(x) \leq M<\infty \hbox{ a.e. } x\in \Omega$, 
\item the set
$\{x\in \Omega|a(x)<0\}$ has positive measure,
\item there exists $c>0$ such that
 $$\int_\Omega (|\nabla u|^2 + a(x) |u|^2 ) \hbox{ } dx
 \geq c\int_\Omega |\nabla u|^2 \hbox{ } dx 
 \hbox{ } \forall u\in H^1_0(\Omega).$$
\end{enumerate}
\noindent In dimension $n=4$ the hypothesis $(2)$
above is replaced by
$$a(x)\leq 0 \hbox{ } a.e. \hbox{ } x\in \Omega.$$

\end{abstract}

\section{Introduction}

In this paper we shall look for the existence of
non--trivial solutions
to the following Yamabe equation:
\begin{equation}\label{yam}-\Delta u+ a(x)u= \mu u|u|^\frac 4{n-2} \hbox{ }  \mu >0,
x\in \Omega \subset {\mathbf R}^n \hbox{ with } n\geq 4\end{equation}
$$u(x)=0 \hbox{ on } \partial \Omega,$$
under suitable assumptions on $a(x)$ that will be specified later.

\vspace{0.1cm}

The main strategy will be is to look at 
the following minimization problem:
\begin{equation}\label{MiN}S_a(\Omega)=\inf_{u\in H^1_0(\Omega)\setminus \{0\}} 
\frac{\int_{\Omega} (|\nabla u|^2 
+ a(x) |u|^2) \hbox{ } dx}{\|u\|_{L^{2^*}}^2},
\end{equation}
where $\Omega\subset {\mathbf R}^n$ is an open set (eventually unbounded),
$\|u\|_{L^{2^*}}^{2^*}=\int_{\Omega} |u|^{2^*} dx$
and $2^*=\frac{2n}{n-2}$. 
Let us underline that the minimization problem 
\eqref{MiN} is non--trivial due to the non--compactness of the Sobolev embedding: 
\begin{equation}\label{crit}H^1_0(\Omega)\subset L^{2^*}(\Omega).
\end{equation}

In fact the problem \eqref{MiN} has been extensively studied in the literature,
starting from the pioneering works \cite{a} and \cite{t},
due to its obvious connection with the equation \eqref{yam}
that plays a fundamental role in Riemannian geometry
(see the very complete book \cite{abook}).

The literature around this problem is too large in order to be exhaustive,
howevere we want to mention at least some of these papers.

\vspace{0.1cm}

In \cite{bn} the problem \eqref{MiN}
has been treated under the following assumptions: $a(x)\in L^\infty(\Omega)$, $\Omega \subset {\mathbf R}^n$
is bounded
and $a(x)\leq -\epsilon<0$ on an open subset of $\Omega$.

\vspace{0.1cm}

In \cite{l} it has been introduced a general approach 
(the concentration--compactness method)
to overcome, in many 
minimization problems, the difficulties 
connected with the lack of compactness in the Sobolev embedding
\eqref{crit}.
In the same paper many applications of the concentration--compactness
method are given, among them let us mention the problem \eqref{MiN} 
that is treated under suitable assumptions on $a(x)$.

Finally we want to mention \cite{te} where the same problem is treated
assuming that $a(x)$ is a function homogeneous
of order $-2$ defined on the whole ${\mathbf R}^n$.

\vspace{0.1cm}

In this article we shall work mainly with functions $a(x)$ belonging to the Lorentz space
$L^{\frac n2, d}(\Omega)$ with $1\leq d <\infty$
(for a definition of $L^{p,q}$ see 
\cite{o} or section \ref{lorentz}) without any further regularity assumption.
Let us point--out that the quadratic form introduced in \eqref{MiN} is meaningful
in general for every $a(x) \in L^{\frac n2, d}(\Omega)$
due to the following inequality:
\begin{equation}\label{intro}
\int_{\Omega} |a(x)| |u|^2 \hbox{ } dx 
\leq C \int_{\Omega}|\nabla u|^2 \hbox{ } dx \hbox{ }
\forall u\in H^1_0(\Omega)
\end{equation}
(see theorem \ref{cpt}) where $C>0$ depends on $a(x)\in L^{\frac n2, d}(\Omega)$.

Notice that if $a(x)\equiv 0$, then the problem
\eqref{MiN} is equivalent to understand whether or not 
the best constant in the critical Sobolev embedding \eqref{crit} is achieved.
By using the concentration--compactness method developed in \cite{l} it is possible to show that
the best constant is achieved when $\Omega\equiv {\mathbf R}^n$.
On the other hand a standard rescaling argument implies that
the best constant is never achieved in the case that $\Omega \neq {\mathbf R}^n$.

\vspace{0.1cm}

As it was mentioned above,
in \cite{bn} the authors have shown that the situation 
changes when a term 
of the type $\int_\Omega a(x) |u|^2 \hbox{ } dx$ 
is added to the energy $\int_\Omega |\nabla u|^2 \hbox{ } dx$, 
provided that $\Omega$ is bounded and 
$a(x)\in L^\infty(\Omega)$ is negative on an open subset of $\Omega$.
 
The main aim of this paper is to show that there exists a minimizer
for \eqref{MiN} when $a(x)$ belongs
to a class more general than the one 
considered in \cite{bn}. Of course for the same class 
of potentials $a(x)$ we can deduce the existence of non--trivial
solutions to \eqref{yam} by using a straightforward 
Lagrange multipliers technique.

Next we state our result in dimension $n\geq5$.

\begin{thm}\label{main}
Let $n\geq 5$ and let $\Omega\subset \mathbf R^n$ an open set
(eventually unbounded). Assume that $a(x)\in L^{\frac n2, d}(\Omega)$
with $d\neq \infty$, satisfies: 
\begin{equation}\label{1a}\hbox{ there exists } 0\leq M<\infty
\hbox{ such that } a(x)< M \hbox{ } a.e. \hbox{ } x\in \Omega;
\end{equation}
\begin{equation}\label{2a} \hbox{ the set } \mathcal N\equiv \{x\in \Omega | a(x)<0\}
\hbox{ has positive measure};
\end{equation}
\begin{equation}\label{3a} \hbox{ there exists } c>0 \hbox{ such that }
\end{equation}
$$\int_{\Omega} (|\nabla u|^2 + a(x)|u|^2)\hbox{ } dx \geq c\int_\Omega |\nabla u|^2 
\hbox{ } dx
\hbox{ } 
\forall u\in H^1_0(\Omega).
$$
Then there exists a function
$v_0\in H^1_0(\Omega)$ such that
$\int_{\Omega} |v_0|^{2^*} dx=1$ and $$\inf_{H^1_0(\Omega)\setminus \{0\}} 
\frac{\int_{\Omega} (|\nabla v|^2 + a(x) |v|^2 )\hbox{ } dx}{\|v\|_{L^{2^*}}^2}
= \int_{\Omega} (|\nabla v_0|^2 + a(x) |v_0|^2 )\hbox{ } dx.
$$
\end{thm}

Let us notice that we do not assume
the continuity of $a(x)$ and
that we allow to the potential $a(x)$ to have a bounded non--negative part.

\vspace{0.1cm}

In dimension $n=4$ we are able to give 
a similar result under an assumption stronger than
\eqref{1a}.
In fact in remark \ref{expint} it is explained where
the proof of theorem \ref{main} fails in dimension $n=4$.
\begin{thm}\label{main4}
Let $\Omega\subset \mathbf R^4$ be an open set
(eventually unbounded). Assume that $a(x)\in L^{2, d}(\Omega)$
with $d\neq \infty$, satisfies \eqref{2a}, \eqref{3a}
as in theorem \ref{main} and 
\begin{equation}\label{4} a(x)\leq 0 \hbox{ } a.e. \hbox{ } x\in \Omega.
\end{equation}
Then there exists a function
$v_0\in H^1_0(\Omega)$ such that
$\int_{\Omega} |v_0|^{2^*} dx=1$ and $$\inf_{H^1_0(\Omega)\setminus \{0\}} 
\frac{\int_{\Omega} (|\nabla v|^2 + a(x) |v|^2 )\hbox{ } dx}{\|v\|_{L^{2^*}}^2}
= \int_{\Omega} (|\nabla v_0|^2 + a(x) |v_0|^2 )\hbox{ } dx.
$$
\end{thm} 
\begin{remark}
Along the proof of theorems \ref{main} and \ref{main4}
it will be clear that the assumption \eqref{3a} is needed only to guarantee
the boundedness in $H^1_0(\Omega)$ of the minimizing sequences in \eqref{MiN}.
On the other hand it is easy to show that
if $a(x)\in L^\frac n2 (\Omega)$, then 
the boundedness of the minimizing sequences
in \eqref{MiN} can be proved by using  the Sobolev embedding
$H^1_0(\Omega)\subset L^{2^*}(\Omega)$ without any 
further assumption. 
This implies that
assumption \eqref{3a} can be removed in theorems \ref{main} 
and \ref{main4} in the case $a(x)\in L^\frac n2(\Omega)$.
On the other hand the coercivity assumption \eqref{3a}
is natural in the literature 
(see for example \cite{bn}, \cite{l}).
\end{remark}

\begin{remark}
Looking at the proof of theorems \ref{main} and \ref{main4} it will be clear that
we prove the following fact:
every minimizing sequence for \eqref{MiN} is compact in $H^1_0(\Omega)$. 
\end{remark}

\begin{remark}\label{previous}
In general theorem \ref{main} cannot be extended to potentials
$a(x)$ belonging  to $L^{\frac n2, \infty}(\Omega)$. For instance it is possible to show that
for every $\mu\in \mathbf R$ and $0<R<\infty$
the value $$H_{\mu, R}=\inf_{H^1_0(|x|<R)\setminus \{0\}} \frac{\int_{|x|<R} 
(|\nabla u|^2 - \mu |x|^{-2}) \hbox{ } dx}{\|u\|_{L^{2^*}}^2}$$
is never achieved.
\end{remark}

\vspace{0.1cm}

We underline that from a technical point of view 
the assumption
$d\neq \infty$, done in the statement of theorems \ref{main} and \ref{main4},
will be relevant in order to prove
the compactness of a Sobolev embedding in suitable weighted spaces.
We think that this result has its own interest and we state it 
separately.

\vspace{0.1cm}

\begin{thm}\label{cpt}
Assume $n\geq 3$ and $\Omega \subset {\mathbf R}^n$. 
Then for every $a(x)\in L^{\frac n2, d}(\Omega)$ with $1\leq d \leq \infty$, we have the following
continuous embedding:
\begin{equation}\label{cont}H^1_0(\Omega) \subset L^2_{|a(x)|}(\Omega)\end{equation}
where
\begin{equation}\label{wei}\|u\|_{L^2_{|a(x)|}(\Omega)}^2=\int_\Omega |a(x)| |u|^2 
\hbox{ } dx.\end{equation}
If moreover $d\neq \infty$, then the embedding is compact.
\end{thm}

Next we fix some notations useful in the sequel.

\vspace{0.1cm}

{\bf Notations.}

For every $1\leq p,q\leq \infty$ we denote by $L^{p,q}(\Omega)$ the usual Lorentz spaces
(see section \ref{lorentz}).

For every $a(x)\in L^{\frac n2, d}(\Omega)$ and $\Omega\subset {\mathbf R}^n$,
we shall denote by $S_a(\Omega)$ the quantity defined in \eqref{MiN}.

We shall make use of the universal constant
\begin{equation}\label{S}
S=\inf_{H^1_0({\mathbf R}^n)
\setminus \{0\}} \frac{\int_{{\mathbf R}^n} |\nabla v|^2 dx}{\|v\|^2_{L^{2^*}}}.
\end{equation}

The norm in the weighted spaces $L^2_{|a(x)|}$ is the one defined in \eqref{wei}.
 
If $A\subset {\mathbf R}^n$ is 
a measurable set then we shall denote by $meas \hbox{ } A$ and $\chi_A$
the measure of $A$ and its characteristic function respectively.

Assume that $X$ is a topological space, then ${\mathcal C}(X)$
denotes the space of continuous and real valued functions 
on $X$.

For every $R>0$ and $x\in {\mathbf R}^n$ we denote by $B_R(x)$
the ball of radius $R$ and centered in $x$.

Given $\alpha\geq 0$ we shall denote by $0(\epsilon^\alpha)$ and $o(\epsilon^\alpha)$
any function of the variable $\epsilon$ such that:
$$\limsup_{\epsilon\rightarrow 0} {|0(\epsilon^\alpha)|}{\epsilon^{-\alpha}}<\infty
\hbox{ and } \lim_{\epsilon\rightarrow 0} {|o(\epsilon^\alpha)|}{\epsilon^{-\alpha}}=0,$$
respectively.

\section{The Lorentz spaces $L^{p,d}(\Omega)$
and proof of theorem \ref{cpt}}\label{lorentz}

In order to introduce the Lorentz spaces
we associate to every measurable 
function its decreasing rearrangement.
Assume that $g:\Omega\rightarrow \mathbf R$ is a measurable function defined on
the measurable set $\Omega\subset {\mathbf R}^n$.
At a first step we associate to the function $g$ its distribution function:
$$m(. , g):(0, \infty]\rightarrow [0, \infty],$$
defined for every $\sigma> 0$ as follows:
$$m(\sigma, g)=meas \hbox{ } \{x\in \Omega| |g(x)|>\sigma\}. $$

It is immediate to show that the distribution function 
defined above is monotonic decreasing.

Once the distribution function $m(\sigma, g)$
has been introduced, we can 
associate to $g$ its decreasing rearrangement function $g^*$:
$$g^*:[0, \infty]\rightarrow [0, \infty],$$
where
$$g^*(t)=inf\{\sigma \in {\mathbf R}^+| m(\sigma, g)<t\}.$$

We can now  define the Lorentz spaces.
\begin{definition}
Assume that $1\leq p< \infty$ and $1\leq d<\infty$, then the measurable function
$g:\Omega \rightarrow \mathbf R$ belongs to
the space $L^{p, d}(\Omega)$ iff
$$\| g\|_{L^{p,d}(\Omega)}^d =\int_0^\infty [g^*(t)]^d t^{\frac dp-1} dt <\infty.$$

\noindent If $1\leq p <\infty$ and $d=\infty$, then
$g\in L^{p, \infty}(\Omega)$ iff
$$\| g\|_{L^{p,\infty}(\Omega)} = \sup_{t>0} g^*(t) t^{\frac 1p} <\infty.$$
\end{definition}

\vspace{0.1cm}

Next we shall describe some
properties satisfied by the functions belonging to the Lorentz spaces
that are important in the sequel (for the proof see \cite{o}).

\begin{prop}
Assume that $1\leq s, s',q,r \leq\infty$ and $1\leq d_1,d_2, d_3 \leq \infty$
are such that:
$$\frac 1s= \frac 1q +\frac 1r \hbox{ and } 
\frac 1s + \frac 1{s'}=1.$$ Then: 

\begin{equation} \label{holder} \int_{\Omega} |fg| dx \leq \|f\|_{L^{s,d_1}(\Omega)} 
\|g\|_{L^{s',d_2}(\Omega)}
\hbox{ provided that } \frac {1}{d_1}+\frac 1{d_2}\geq 1;
\end{equation}
\begin{equation}\label{Holder} \|f g\|_{L^{s, d_1}(\Omega)}\leq s' \|f\|_{L^{q,d_2}(\Omega)} 
\|g\|_{L^{r,d_3}(\Omega)} 
\hbox{ provided that } \frac 1{d_2} + \frac 1{d_3}\geq \frac 1{d_1}.
\end{equation}
\end{prop}

\vspace{0.1cm}

Next result is a well--known improved version  
of the classical Sobolev embedding (see \cite{m} and \cite{t}).

\begin{prop}\label{lorentzSob}
For every $n\geq 3$ and for every
open set $\Omega \subset {\mathbf R}^n$, 
there exists a real constant $C=C(\Omega)>0$
such that:
\begin{equation}\label{sobolev}
\| f\|_{L^{2^*,2}(\Omega)}\leq C \|\nabla f\|_{L^2(\Omega)}
\hbox{  } \forall f\in H^1_0(\Omega).
\end{equation}
\end{prop}

\begin{remark}
We want to underline that \eqref{sobolev} represents an improved 
version of the standard Sobolev embedding
due to the following inclusions:
$$L^{p,d_1}\subset L^{p,p}=L^p\subset L^{p, d_2}$$
where $1\leq d_1\leq p \leq d_2 \leq \infty$.
\end{remark}

We are now able to prove theorem \ref{cpt}.

\noindent {\bf Proof of thorem \ref{cpt}.}
First we prove the continuity of the embedding \eqref{cont}.
Notice that since $L^{\frac n2, d}(\Omega)\subset L^{\frac n2, \infty}(\Omega)$
for every $1\leq d \leq \infty$, we can assume $d=\infty$.
By combining \eqref{holder}, \eqref{Holder} and \eqref{sobolev}
we get:
\begin{equation}\label{LuC}\int_{\Omega} |a(x)| |u|^2 \hbox{ } dx\leq 
\|a(x)\|_{L^{\frac n2, \infty}(\Omega)} \|u^2\|_{L^{\frac n{n-2}, 1}(\Omega)}
\end{equation}
$$\leq \frac n2 \|a(x)\|_{L^{\frac n2, \infty}(\Omega)}
 \|u\|_{L^{2^*, 2}(\Omega)}^2
\leq C \|a(x)\|_{L^{\frac{n}{2}, \infty}(\Omega)}
\|u\|_{H^1_0(\Omega)}^2.
$$

\vspace{0.1cm}

Next we prove the compactness of the embedding 
\eqref{cont} when $d\neq \infty$.

\vspace{0.1cm}

Let
$\{u_k\}_{k\in \mathbf N}$ be a sequence bounded in 
$H^1_0(\Omega)$. We can assume that up to a subsequence
there exists a function $u_0\in H^1_0(\Omega)$ such that:
$$u_k\rightharpoonup u_0\hbox { in } H^1_0(\Omega).$$ 

We shall show that up to a subsequence we have:
\begin{equation}\label{FIN}
\lim_{k\rightarrow \infty}
\int_{\Omega} |a(x)| |u_k|^2 \hbox{ } dx
=\int_\Omega |a(x)| |u_0|^2 \hbox{ } dx,\end{equation}
and it will complete the proof.
To show \eqref{FIN} let us first notice the following property:
\begin{align}\label{bound} 
&\hbox{{\em  for every bounded  open set} } K \hbox{{\em such that }}
 K\subset  \Omega, 
\hbox{{\em there exists }} C=C( K)>0 
\\ \nonumber &\hbox{ {\em such that }} \|u_k\|_{H^1( K)}< C.
\end{align}

In fact it is sufficient to show that the $L^2$--norm 
of the functions $\{u_k\}_{k\in \mathbf N}$ are bounded on every bounded set,
since the boundedness of the $L^2$--norm of the gradients
comes from the assumption. 

Due to 
the H\"older inequality and to the Sobolev embedding we have: 
$$\|u_k\|_{L^2(K)}
\leq  |K|^\frac 1n 
\|u_k\|_{L^\frac{2n}{n-2}(\Omega)}\leq C \| \nabla u_k\|_{L^2(\Omega)}$$
where $C>0$ is a suitable constant
that depends on $K$,
then 
\eqref{bound} holds.

For every  $i\in \mathbf N$ 
we split the domain $\Omega$ as follows: 
$$\Omega=\Omega^{i}_1 \cup \Omega^{i}_2,$$ 
where the splitting is the one described in 
proposition \ref{propr} and corresponding to 
$\epsilon=\frac 1i$ (see the Appendix).

Following the proof of \eqref{LuC}
one can deduce that there exists a constant $C>0$ that depends only on $\Omega$ and such that:
\begin{align}\label{infi}
&\int_{\Omega^{i}_2} 
|a(x)| |u_k|^2 \hbox{ } dx \leq 
C \|\chi_{\Omega^{i}_2}|a(x)|\|_{L^{\frac{n}{2}, d}(\Omega)}
\|\nabla u_k\|_{L^{2}(\Omega)}^2 < \frac C i \hbox{  } \forall i, k \in \mathbf N,
\end{align} 
where we used the boundedness 
of the sequence $\{u_k\}_{k\in \mathbf N}$ in $H^1_0(\Omega)$, 
and the properties of $\Omega^{i}_2$
described in proposition \ref{propr}.

Recall also that $\{\Omega^{i}_1\}_{i\in \mathbf N}$
is a sequence of bounded domains.
We can then combine \eqref{bound} with
the compactness of the Sobolev embedding
on the bounded domain in order to deduce that:
$$\|u_k - u_0\|_{L^2(\Omega^{i}_1)}\rightarrow 0 \hbox{ as } k\rightarrow \infty,$$
where $i\in \mathbf N$ is a fixed number.
Due to propositon \ref{propr} we have that $|a(x)|$ is bounded on 
$\Omega^{i}_1$, then the previous inequality implies that
for every $i\in \mathbf N$ there exists $k(i)\in \mathbf N$ such that:
 
\begin{equation}\label{comp}\left |\int_{\Omega^{i}_1} |a(x)||u_{k(i)}|^2 \hbox{ } dx
- \int_{\Omega^{i}_1} |a(x)| |u_0|^2 \hbox{ } dx \right | 
<  \frac 1 i.
\end{equation}

It is easy to show that in fact we can choose $k(i)$ in such a way that
$k(i)<k(i+1)$.
By combining \eqref{infi}
with \eqref{comp},
and by using a diagonalizing argument, we can conclude that 
up to a subsequence we have:
$$\lim_{k\rightarrow \infty}\int_{\Omega} 
|a(x)| |u_k|^2 \hbox{ } dx=\int_{\Omega} |a(x)||u_0|^2 \hbox{ } dx.$$

\hfill$\Box$

\section{A general approach to the problem \eqref{MiN}}\label{brni}

Let $v_n\in H^1_0(\Omega)$ be a sequence such that:
$$\int_\Omega |v_n|^{2^*} \hbox{ } dx=1$$
and
\begin{equation}\label{bet}
\lim_{n\rightarrow \infty}\int_{\Omega} (|\nabla v_n|^2 + a(x) |v_n|^2 )
\hbox{ } dx=S_a(\Omega)
\end{equation}
where $a(x)$ satisfies the assumptions done in theorems \ref{main} and \ref{main4}.
Notice that due to assumption \eqref{3a}
we can deduce that
$\{v_n\}_{n\in \mathbf N}$ is bounded in $H^1_0(\Omega)$.

\vspace{0.1cm}

Moreover the weak--compactness of bounded sequences in $H^1_0(\Omega)$
and the compactness of the embedding
given in theorem \ref{cpt}, imply the existence of
$v_0\in H^1_0(\Omega)$ such that up to a subsequence we have:

\vspace{0.1cm}

\begin{enumerate}
\item $v_n\rightharpoonup v_0$ in $H^1_0(\Omega)$;
\item $\lim_{n\rightarrow \infty}\int_{\Omega} a(x) |v_n|^2 
\hbox{ } dx= 
\int_\Omega a(x)|v_0|^2 \hbox{ } dx.$
\end{enumerate}

Notice that by combining $(1)$ and $(2)$ with \eqref{bet} one deduce that
\begin{equation}\label{semic}
\int_{\Omega} (|\nabla v_0|^2 + a(x)|v_0|^2) \hbox{ } dx\leq S_a(\Omega).
\end{equation}
\vspace{0.1cm}

On the other hand, following the same argument as in Br\'ezis and Nirenberg
(see also \cite{a} and \cite{l})
and recalling $(1)$ and $(2)$ above, one can deduce the following implication:
\begin{equation}\label{basic}
\hbox{ if }S_a(\Omega)<S
\hbox { then } \lim_{n\rightarrow \infty}\|v_n-v_0\|_{L^{2^*}(\Omega)}=0
\end{equation}
(recall that $S$ is defined in \eqref{S}). In
particular if $S_a(\Omega)<S$ then  $\int_{\Omega}|v_0|^{2^*}=1$, and in turn
it can be combined with
\eqref{semic} to deduce  
that the value $S_a(\Omega)$ is achieved in $H^1_0(\Omega)\setminus \{0\}$
when $S_a(\Omega)<S$ .

The main purpose in next sections will be to prove that $S_a(\Omega)<S$ under 
the assumptions done on $a(x)$ in theorems \ref{main} and \ref{main4}.

\vspace{0.1cm}

Next we recall a basic fact proved in \cite{bn} 
that will be the starting point in the proof
(at least in the case $n>4$) of the inequality $S_a(\Omega)<S$.

\vspace{0.1cm}

Assume that $n\in \mathbf N$ is fixed. 
We shall denote by 
$u_\epsilon(x)$ the following family of rescaled functions:
\begin{equation}\label{resc}
u_\epsilon(x)=\frac{[n(n-2) \epsilon^2]^\frac{n-2}{4}}{(\epsilon^2 + |x|^2)^\frac{n-2}{2}}
\hbox{ }\forall x \in {\mathbf R}^n, \epsilon>0\end{equation}
and for every $x_0\in {\mathbf R}^n$
$$u_{\epsilon, x_0}=u_\epsilon(x-x_0).$$
Let us recall that the functions $u_\epsilon$ defined above 
realize the best constant in the critical Sobolev embedding (see \cite{tal}).
In fact 
it is possible to prove that
\begin{equation}\label{dirichlet}
\int_{{\mathbf R}^n} |\nabla u_\epsilon|^2 \hbox{ } dx=S^\frac n2
\end{equation}
\begin{equation}\label{mass2}
\int_{{\mathbf R}^n} |u_{\epsilon}|^{2^*} \hbox{ } dx = S^\frac n2 
\end{equation}
for every $n\geq 3$.

Let us fix also a cut--off function $\eta\in C^\infty_{0}(|x|<2)$ such that
$0\leq \eta \leq 1$ and $\eta\equiv 1$ in $\{|x|<1\}$.
For every $\mu>0, x_0 \in {\mathbf R}^n$
we introduce the function 
$$\eta_{\mu, x_0}(x)=\eta \left( \frac{x-x_0}{\mu}\right).$$

Next we state a basic proposition whose proof can be found in \cite{bn} (see also 
\cite{a} and \cite{st}).

\begin{prop}\label{uepsilon}
Let $\lambda<0$ be a fixed number and $\Omega \subset {\mathbf R}^n$ with $n\geq 5$ be
an open set.
For any
$\mu>0$, $x_0\in \Omega$
there exists $c=c(n)>0$  
such that the following estimate holds:
\begin{equation}\label{nap} \int_{{\mathbf R}^n}(|\nabla (u_{\epsilon, x_0}
\eta_{\mu, x_0}) |^2 + \lambda |u_{\epsilon, x_0} \eta_{\mu, x_0}|^2 ) \hbox{ } dx
\leq S^\frac n2 + c \lambda \epsilon^2+0(\epsilon^{n-2}) 
\end{equation}
(here $0(\epsilon^{n-2})$ depends on $\mu>0$).\\
\noindent Moreover for every $n\geq 4$ and for every $\mu>0$ we have
\begin{equation}\label{mass}
\|u_{\epsilon, x_0} \eta_{\mu, x_0}\|_{L^{2^*}}^{2^*}=S^\frac n2 + 0(\epsilon^n)
\end{equation}
(here $0(\epsilon^n)$ depends on $\mu$).
\end{prop}

\begin{remark}
If $n=4$ then  the following asymptotic behaviour
is given in \cite{bn}:  
\begin{equation}\label{nap4}\int_{{\mathbf R}^4} (|\nabla (u_{\epsilon, x_0}\eta_{\mu, x_0}) 
|^2 + \lambda |u_{\epsilon, x_0}
\eta_{\mu, x_0}|^2 ) \hbox{ } dx
\leq  S^2 + c \lambda \epsilon^2 |\ln \epsilon|+0(\epsilon^{2}). 
\end{equation}
In fact along the proof of theorem \ref{main4} (more precisely in the proof of lemma
\ref{lem4}) we shall need a slightly refined version of this estimate.
\end{remark}

\begin{remark}\label{expint}
Notice that 
$u_\epsilon(x)\in L^2(\mathbf R^n)$ for every $n>4$ and
$u_\epsilon (x)\notin L^2(\mathbf R^4)$ in the case $n=4$. 
In next sections it will be clear that this is the main 
reason why the dimension $n=4$ will be treated in a different way
compared with the dimensions $n>4$. 
\end{remark}

\section{Proof of theorem \ref{main}}

In this section the functions $u_{\epsilon, x_0}$
and $\eta_{\mu, x_0}$ are the ones introduced in section \ref{brni}. 

Let us recall also that 
in order to prove theorem \ref{main} it 
is sufficient to prove the following lemma (see section \ref{brni}).
\begin{lem}\label{n5}
Assume $n\geq 5$ and $\Omega\subset {\mathbf R}^n$. If $a(x)$ 
satisfies the assumptions of theorem \ref{main}
then $S_a(\Omega)<S$.
\end{lem}
\noindent {\bf Proof.}
Since now on we assume that a representative of the function $a(x)$
has been fixed
and we shall not consider $a(x)$ as a class of functions that are
equivalent modulo zero measure sets.
This will allow us to consider the pointwise value $a(x)$ for every fixed $x\in \Omega$.

First of all we notice that we can assume $a(x)\in L^\infty(\Omega)\cap L^{\frac n2, d}(\Omega)$.

In fact it is easy to show that assumption \eqref{2a} implies that
there exists $N_0\in \mathbf N$ such that 
$$meas \hbox{ } \{ x \in \Omega| -N_0<a(x)<0\}>0.$$

In particular the potential $$\tilde a(x)=Max\{a(x), -N_0\}$$
satisfies the assumptions of theorem \ref{main} and moreover
$\tilde a(x)\in L^\infty(\Omega)\cap L^{\frac n2, d}(\Omega)$ 
(notice that this is stronger than \eqref{1a}).
By combining this fact with the following trivial inequality:  
$$S_{a}(\Omega) \leq S_{\tilde a}(\Omega),$$
one deduce that it is not restrictive 
to assume $a(x)\in L^{\infty}(\Omega) \cap L^{\frac n2, d}(\Omega)$.
In particular we can assume that $a(x)\in L^{\frac n2}_{loc}(\Omega)$.

We are then in position to use the Lebesgue derivation theorem
in order to deduce that (see \cite{eg} for a proof):
\begin{equation*}
\lim_{\epsilon \rightarrow 0} \epsilon^{-n} \int _{B_\epsilon(\bar x)}
|a(x) - a(\bar x)|^\frac n2 \hbox{ } dx=0 \hbox{ }  a.e. \hbox{ } \bar x\in \Omega,
\end{equation*}
that due to assumption \eqref{2a} implies
the existence of $x_0\in \Omega$ such that
\begin{equation}\label{leb}
\lim_{\epsilon \rightarrow 0} \epsilon^{-n} \int _{B_\epsilon(x_0)}
|a(x) - a(x_0)|^\frac n2 \hbox{ } dx=0 \hbox{ and } -\infty<a(x_0) <0.
\end{equation}

Due again to the the definition of $S_a(\Omega)$ it is easy to verify that
the following inequality holds:
$$S_{a}(\Omega)\leq S_{Max\{a(x), a(x_0)\}}(\Omega),$$
and it implies clearly that it not restrictive to assume that:
\begin{equation}\label{rest}a(x_0)\leq  a(x)\leq M \hbox{ } a.e. 
\hbox{ } x \in \Omega
\end{equation}
(here we have used \eqref{1a} in the r.h.s. inequality).

Next we notice that the following identity holds trivially:
\begin{equation}
\label{uno}\int_{\Omega} (|\nabla (u_{\epsilon, x_0} \eta_{\mu, x_0})|^2 + 
a(x)|u_{\epsilon, x_0}\eta_{\mu, x_0} |^2) \hbox{ } dx
\end{equation}
$$
=\int_{\Omega} (|\nabla (u_{\epsilon,x_0} \eta_{\mu, x_0})|^2 + 
a(x_0) |u_{\epsilon, x_0}\eta_{\mu, x_0}|^2) \hbox{ } dx +
\int_{\Omega}(-a(x_0) + a(x))|u_{\epsilon, x_0} \eta_{\mu, x_0}|^2 \hbox{ } dx
$$
where $\mu>0$ is choosen small enough 
in such a way that $u_{\epsilon, x_0}\eta_{\mu, x_0}\in H^1_0(\Omega)$.

\vspace{0.1cm}

By using the H\"older inequality we get:
\begin{equation}\label{due}
\left |\int_{\Omega}(- a(x_0) + a(x))|u_{\epsilon, x_0} \eta_{\mu, x_0}|^2 
\hbox{ } dx\right|
\end{equation}$$
\leq \int_{\Omega\cap B_{\epsilon R}(x_0)}|a(x_0) - a(x)||u_{\epsilon, x_0}\eta_{\mu, x_0}|^2  
\hbox{ } dx 
$$$$+ \int_{\Omega \cap \{|x-x_0|>\epsilon R\}} |a(x_0) - a(x)||u_{\epsilon, x_0} \eta_{\mu, x_0}|^2 
\hbox{ } dx $$
$$\leq \|u_\epsilon\|_{L^{2^*}({\mathbf R}^n)}^2\left (\int_{\Omega \cap B_{\epsilon R}(x_0)}
|a(x_0) - a(x)|^\frac n2 \hbox{ } dx\right )^\frac 2n 
$$$$+ 2 Max \hbox{ } \{|a(x_0)|, M\} \int_{|x|>\epsilon R} |u_\epsilon|^2 \hbox{ } dx$$
where we have used \eqref{rest}
and $R>0$ is a 
real number that we are going to fix.

In fact we choose $R>0$ large enough in such a way that the following condition
holds:
\begin{equation}\label{tre}\int_{|x|>\epsilon R} |u_\epsilon|^2 
\hbox{ } dx
=
\frac{[n(n-2)]^\frac{n-2}{2}}{\epsilon^{n-2}} 
\int_{|x|>\epsilon R}
\frac{dx}{ [1 + \epsilon^{-2}|x|^2]^{n-2}}
\end{equation}$$= [n(n-2)]^\frac{n-2}{2}\epsilon^{2} 
\int_{|x|>R}
\frac{dx}{ [1 + |x|^2]^{n-2}}<\frac {c \hbox{ } |a(x_0)|}{4 Max \hbox{ }\{|a(x_0)|, M\}} \epsilon^2
$$
where $c>0$ is the same constant that appears in \eqref{nap}.

On the other hand due to \eqref{leb} we deduce that
\begin{equation}
\label{quattro}\lim_{\epsilon\rightarrow 0} R^{-2}\epsilon^{-2}\left(\int_{B_{\epsilon R}(x_0)}
|a(x_0) - a(x)|^\frac n2 \hbox{ } dx\right )^\frac 2n
=0.\end{equation}

By combining \eqref{uno}, \eqref{due}, \eqref{tre}, \eqref{quattro}
with \eqref{nap} we get:
$$\int_{\Omega} (|\nabla (u_{\epsilon, x_0} \eta_{\mu, x_0})|^2 + a(x)|u_{\epsilon, x_0} \eta_{\mu, x_0}|^2) 
\hbox{ } dx
$$$$\leq S^\frac n2 + c \hbox{ } a(x_0) \epsilon^2 + \frac c2 \hbox{ } |a(x_0)| \epsilon^2
+ S^\frac{n-2}2 R^2 \epsilon^2 o(1) + 0(\epsilon^{n-2}) $$
$$=S^\frac n2 +\frac c2 \hbox{ }a(x_0) \epsilon^2 + o(\epsilon^2).$$

By using now \eqref{mass} we deduce that for $\epsilon>0$ small enough
the following chain of inequalities holds:
$$S_a(\Omega)\leq \frac{\int_{\Omega} (|\nabla (u_{\epsilon, x_0} \eta_{\mu, x_0})|^2 
+ a(x)|u_{\epsilon, x_0} \eta_{\mu, x_0} |^2) \hbox{ } dx}
{\|u_{\epsilon,x_0} \eta_{\mu, x_0}\|_{L^{2^*}(\Omega)}^2}$$
$$\leq \frac{S^\frac n2 +\frac c2 
\hbox{ } a(x_0) \epsilon^2 + o(\epsilon^2)}
{(S^\frac n2 + 0(\epsilon^n))^\frac 2{2^*}}<S,$$
where at the last step we have used that $ c \hbox{ } a(x_0)<0$.

\hfill$\Box$

\section{Proof of theorem \ref{main4}}

Notice that for $n=4$ the functions $u_\epsilon(x)$
(that have been introduced in section \ref{brni}) become:
$$u_\epsilon(x)=\frac{\sqrt 8 \hbox{ } \epsilon }
{\epsilon^2 + |x|^2} \hbox{ } \forall x\in {\mathbf R}^4\hbox{ } \forall \epsilon>0.$$

We shall also need
$$u_{\epsilon, x_0}=u_\epsilon(x-x_0)$$
and 
$$\eta_{\mu, x_0}=\eta\left(\frac{x-x_0}{\mu} \right)$$
where $x_0\in {\mathbf R}^4$, $\mu>0$ and $\eta$ is a cut--off function
belonging to $C^\infty_0(|x|<2)$ such that
$0\leq \eta \leq 1$ and $\eta\equiv 1$ in $\{|x|<1\}$.
 
Let us recall also that in dimension $n=4$ 
the identity
\eqref{mass} becomes:
\begin{equation}\label{mass4}
\int_{{\mathbf R}^4}|u_{\epsilon, x_0}\eta_{\mu, x_0}|^4 \hbox{ } dx
=S^2 + 0(\epsilon^4)
\end{equation}
where $0(\epsilon^4)$ depends on $\mu>0$.

\vspace{0.1cm}

Next we prove a lemma that is sufficient
to conclude the proof of theorem \ref{main4}
(see section \ref{brni}).

\begin{lem}\label{lem4}
Assume that $\Omega \subset \mathbf R^4$ is an open set and $a(x)$ satisfies the 
assumptions of theorem \ref{main4}
then $S_a(\Omega)<S$.
\end{lem}

\noindent {\bf Proof.}
As in the proof of lemma \ref{n5}
we assume that a representative of the function $a(x)$
has been fixed.
This will allow us to consider the pointwise value $a(x)$ for every fixed $x\in \Omega$.

Arguing as in the first part of the proof of lemma \ref{n5} we can 
assume the existence of $x_0\in \Omega$ such that:
\begin{equation}\label{leb4}
\lim_{\epsilon \rightarrow 0} \epsilon^{-4} \int _{B_\epsilon(x_0)}
|a(x) - a(x_0)|^2 \hbox{ } dx=0 \hbox{ with } -\infty<a(x_0) <0 
\end{equation}
and
\begin{equation}\label{nores}
a(x_0)\leq a(x)\leq 0 \hbox{ } a.e. \hbox{ } x\in \Omega
\end{equation}
(here we have used \eqref{4} in the r.h.s inequality).

Since now on we fix $\mu>0$ such that
$$supp \hbox{ } \eta_{\mu, x_0}\subset \Omega$$ 
(in fact this condition is sufficient to deduce that
$u_{\epsilon, x_0}\eta_{\mu, x_0}\in H^1_0(\Omega)$).

Let us write
the following trivial identity:
\begin{equation}
\label{uno4}\int_{\Omega} (|\nabla (u_{\epsilon,x_0} \eta_{\mu,x_0})|^2 
+ a(x)|u_{\epsilon,x_0}\eta_{\mu,x_0} |^2) \hbox{ } dx
=I_\epsilon^{\mu}+II_\epsilon^{\mu} \end{equation}
where
$$
I_\epsilon^{\mu}=\int_{\Omega} (|\nabla (u_{\epsilon, x_0}\eta_{\mu, x_0})|^2 + 
a(x_0) |u_{\epsilon, x_0}\eta_{\mu, x_0}|^2) \hbox{ } dx 
$$$$
II_\epsilon^{\mu}
=\int_{\Omega}(-a(x_0) + a(x))|u_{\epsilon, x_0}\eta_{\mu, x_0}|^2 \hbox{ } dx.
$$

\vspace{0.1cm}

{\em Estimate for $I_\epsilon^{\mu}$}

\vspace{0.1cm}

Notice that we have:$$\int_{\Omega} |\nabla (u_{\epsilon, x_0}\eta_{\mu, x_0})|^2 \hbox{ } dx=
\int_\Omega |\eta_{\mu, x_0} \nabla u_{\epsilon, x_0} + u_{\epsilon, x_0}
\nabla \eta_{\mu, x_0}|^2 \hbox{ } dx$$
$$\leq \int_\Omega |\eta_{\mu, x_0} \nabla u_{\epsilon, x_0}|^2 \hbox{ } dx
+ \int_\Omega |u_{\epsilon, x_0}
\nabla \eta_{\mu, x_0}|^2 \hbox{ } dx $$$$+ 2 \int_\Omega \eta_{\mu, x_0}
|\nabla \eta_{\mu, x_0}||\nabla u_{\epsilon, x_0}|
|u_{\epsilon, x_0}| \hbox{ } dx$$
$$\leq 32 \epsilon^2 \int_{|x|< 2 \mu}
\frac{|x|^2 dx}{(\epsilon^2 + |x|^2)^4} 
+ \frac{8 \epsilon^2}{\mu^2} \|\nabla \eta\|_{L^\infty({\mathbf R}^4)}^2
\int_{\mu < |x|< 2 \mu } \frac{dx}{(\epsilon^2 + |x|^2)^2}$$
$$+ \frac{32 \hbox{ } \epsilon^2}{\mu} \|\nabla \eta\|_{L^\infty({\mathbf R}^4)}
\int_{\mu < |x|<2 \mu } \frac{ |x| dx}{(\epsilon^2 + |x|^2)^3}$$
and then with elementary computations 
\begin{equation}\label{stup27}
\int_{\Omega} |\nabla (u_{\epsilon, x_0}\eta_{\mu, x_0})|^2 \hbox{ } dx
\leq \int_{{\mathbf R}^4} |\nabla u_\epsilon|^2 \hbox{ } dx -
32 \int_{|x|>2 \mu} 
\frac{\epsilon^2 |x|^2 dx}{(\epsilon^2 + |x|^2)^4} 
\end{equation}
$$+ \left ( 32 \omega_3 \|\nabla \eta\|_{L^\infty({\mathbf R}^n)}^2  
+ 256 \omega_3 \|\nabla \eta\|_{L^\infty({\mathbf R}^n)} \right )\frac{\epsilon^2}{\mu^2},
$$
where $\omega_3$ denotes the Haussdorf measure of the sphere 
${\mathcal S}^3$.
\vspace{0.1cm}

Notice that due to \eqref{dirichlet} the inequality \eqref{stup27}
implies 
\begin{equation}\label{stup}
\int_{\Omega} |\nabla (u_{\epsilon, x_0}\eta_{\mu, x_0})|^2 \hbox{ } dx
\leq S^2  
+ C\frac{\epsilon^2}{\mu^2},
\end{equation}
where $C>0$ is an universal constant.

On the other hand
\begin{equation}\label{stup2}\int_\Omega |u_{\epsilon,x_0}\eta_{\mu,x_0} |^2 
\hbox{ } dx\geq 
8 \epsilon^2 \int_{|x|<\epsilon R} \frac{dx}{(\epsilon^2+|x|^2)^2}
\end{equation}
$$+ 8\epsilon^2 \int_{\epsilon R<|x|<\mu} \frac{dx}{(\epsilon^2+|x|^2)^2}$$
where $R>0$ is a number that we shall fix later
and $\epsilon<\frac \mu R$
(recall that $\mu$ has been fixed above).

\vspace{0.1cm}

By combining \eqref{stup} with \eqref{stup2} and recalling that $a(x_0)<0$
we get:
\begin{equation}\label{stup23}
I_\epsilon^{\mu} \leq 
S^2  
+ C \frac{\epsilon^2}{\mu^2}\end{equation}
$$
+ 8 a(x_0) \epsilon^2 \int_{|x|<\epsilon R} \frac{dx}{(\epsilon^2+|x|^2)^2}
+ 8 a(x_0) \epsilon^2 \int_{\epsilon R<|x|<\mu} \frac{dx}{(\epsilon^2+|x|^2)^2}$$
where $C>0$ is an universal constant.

\vspace{0.1cm}

{\em Estimate for $II_\epsilon^{\mu}$}

\vspace{0.1cm}

The H\"older inequality implies:
\begin{equation}\label{due4}
\left |\int_{\Omega}(-a(x_0) + a(x))|u_{\epsilon, x_0} \eta_{\mu, x_0}|^2 \hbox{ } dx\right|
\end{equation}
$$\leq \|u_\epsilon\|_{L^{4}({\mathbf R}^4)}^2\left (\int_{B_{\epsilon R}(x_0)}
|a(x_0) - a(x)|^2 \hbox{ } dx\right )^\frac 12 
+ |a(x_0)| \int_{\epsilon R <|x|< 2 \mu } 
|u_\epsilon|^2 \hbox{ } dx,$$
where we have used \eqref{nores} to deduce $|a(x)-a(x_0)|\leq |a(x_0)|$,
while $R>0$ is a constant that we shall fix later.

\vspace{0.1cm}

Notice that \eqref{leb4} implies
\begin{equation}
\label{quattro4}\lim_{\epsilon\rightarrow 0} 
{(\epsilon R)}^{-2}\left(\int_{B_{\epsilon R}(x_0)}
|a(x_0) - a(x)|^2 dx\right )^\frac 12
=0,\end{equation}
while the definition of $u_\epsilon$ gives
\begin{equation}\label{tre4}\int_{\epsilon R<|x|<2 \mu} |u_\epsilon|^2 dx
= 8\epsilon^{2} 
\int_{\epsilon R <|x|< 2 \mu }
\frac{dx}{(\epsilon^2 + |x|^2)^{2}}.
\end{equation}

By combining \eqref{mass2}, \eqref{due4}, \eqref{quattro4} and \eqref{tre4} we get:
\begin{equation}\label{partial2}
|II_\epsilon^{\mu}|\leq S R^2 \epsilon^2 \hbox{ } o(1)
+ 8 \epsilon^2  \hbox{  } |a(x_0)|
\int_{\epsilon R <|x|<2 \mu }
\frac{dx}{(\epsilon^2 + |x|^2)^{2}}.\end{equation}

Due to \eqref{uno4}, \eqref{stup23} and \eqref{partial2} we finally
get:
$$\int_{\Omega} (|\nabla (u_{\epsilon,x_0} \eta_{\mu,x_0})|^2 
+ a(x)|u_{\epsilon,x_0}\eta_{\mu,x_0} |^2) \hbox{ } dx\leq S^2 
+ C \frac{\epsilon^2}{\mu^2}$$$$
+ 8 a(x_0) \epsilon^2 \int_{|x|<\epsilon R} \frac{dx}{(\epsilon^2+|x|^2)^2}
+ 8 |a(x_0)| \epsilon^2 \int_{\mu<|x|<2 \mu} 
\frac{dx}{(\epsilon^2+|x|^2)^2}
+ R^2 o(\epsilon^2)$$
$$\leq 
 S^2 + C \frac{\epsilon^2}{\mu^2}
+ 8 a(x_0) \epsilon^2 \int_{|x|<R} \frac{dx}{(1+|x|^2)^2}
+ 8 \omega_3 
|a(x_0)| \epsilon^2 \ln 2 + R^2 o(\epsilon^2)$$
where $C>0$ is an universal constant, $\omega_3$ is the measure of the sphere 
${\mathcal S}^3$
and $R>0$ is a number to be fixed later.
 
Then we have proved the following estimate:
\begin{equation}\label{,mat}\int_{\Omega} (|\nabla (u_{\epsilon,x_0} \eta_{\mu,x_0})|^2 
+ a(x)|u_{\epsilon,x_0}\eta_{\mu,x_0} |^2) \hbox{ } dx\end{equation}
$$\leq S^2 + 
\left( 8\phi(R) a(x_0)+\frac{C}{\mu^2}+8\omega_3 |a(x_0)| \ln 2 \right)\epsilon^2 + R^2 o(\epsilon^2),
$$
where $$\phi(R)=\int_{|x|<R} \frac{dx}{(1+|x|^2)^2}$$
and hence
\begin{equation}\label{gran}
\lim_{R\rightarrow \infty} \phi(R)=\infty.
\end{equation} 

Due to 
\eqref{mass4} and \eqref{,mat}
we deduce that for $\epsilon>0$ small enough we get:
\begin{equation}\label{test}
S_a(\Omega)\leq \frac{\int_{\Omega} (|\nabla (u_{\epsilon, x_0} \eta_{\mu, x_0})|^2 
+ a(x)|u_{\epsilon, x_0} \eta_{\mu, x_0} |^2) \hbox{ } dx}
{\|u_{\epsilon,x_0} \eta_{\mu, x_0}\|_{L^{2^*}(\Omega)}^2}\end{equation}
$$\leq \frac{S^2 + \left( 8\phi(R) a(x_0)+\frac{C}{\mu^2}+8\omega_3 |a(x_0)| \ln 2 \right)
\epsilon^2 + R^2 o(\epsilon^2)}
{(S^2 + 0(\epsilon^4))^\frac 1{2}}<S $$
where we have used at the last step that $ a(x_0) <0$
and we are assuming that  $R>0$ is large enough
in order to guarantee that
$8\phi(R)a(x_0)+\frac C{\mu^2}+ 8\omega_3 |a(x_0)| \ln 2<0$ 
(note that it is possible
due to \eqref{gran}).   

\hfill$\Box$

\section{Appendix}

In order to make this article self--contained we give 
the proof of a proposition contained in \cite{v}. We recall
also that next result has been fundamental along the proof
of theorem \ref{cpt}.
 
\begin{prop}\label{propr}
Let $n\geq 1$ and $\Omega \subset {\mathbf R}^n$ an open set.
Assume that $a(x)\in L^{p, d}(\Omega)$ for $1\leq p, d<\infty$.
Then for any $\epsilon>0$ there exist two measurable sets
$\Omega^{\epsilon}_1, \Omega^{\epsilon}_2$ such that:
$$\Omega^{\epsilon}_1 \cup \Omega^{\epsilon}_2=\Omega,
\Omega^{\epsilon}_1 \cap \Omega^{\epsilon}_2
=\emptyset, \Omega^{\epsilon}_1 \hbox{ is bounded }$$  
and  $$a(x)\chi_{\Omega^{\epsilon}_1}  \in L^\infty(\Omega), \hbox{ } 
\|a(x) \chi_{\Omega^{\epsilon}_2}\|_{L^{p, d}(\Omega)}<\epsilon.$$
\end{prop}

In next lemma the function $f^*$ associated to a function $f$ is the one defined in section
\ref{lorentz}.

\begin{lem}\label{monotono}
Assume that $f_k:\Omega \rightarrow \mathbf R$
is a sequence of functions such that:
$$f_k(x)\geq 0 \hbox { a.e. } x\in \Omega, \hbox{ }f_1\in L^{p, d}(\Omega) 
\hbox{ for suitable } 1\leq p <\infty, 1\leq d<\infty,
$$\begin{equation}\label{MONo}
0\leq  f_{k+1}(x)\leq f_{k}(x) 
\hbox{ a.e. } x\in \Omega \hbox{  } \forall k\in \mathbf N,
\end{equation}
\begin{equation}\label{limpun} \hbox { and }
\lim_{k\rightarrow \infty} f_k(x)=0 \hbox{ a.e. } x\in \Omega.
\end{equation}
Then 
\begin{equation}
\label{cuppe}f_{k+1}^*(t)\leq f_{k}^*(t)
\hbox{ } \forall t\in {\mathbf R}^+\hbox{ } \forall k\in \mathbf N
\end{equation}
and
\begin{equation}\label{null} \lim_{k\rightarrow \infty} 
f_k^*(t)=0 \hbox{  }  \forall t\in {\mathbf R}^+.
\end{equation}
\end{lem}

\noindent {\bf Proof.}
The assumption $0\leq f_{k+1}\leq f_k$ implies that:
$$m(\sigma, f_{k+1}) \leq m(\sigma, f_{k}),$$
where $m(\sigma, g)$ is defined as
in section \ref{lorentz}
for every measurable function $g$. 
Due to this inequality and to the definition of
$f_k^*$ it is easy to deduce
that $f_{k+1}^*(t)\leq f_{k}^*(t)$ and hence \eqref{cuppe}
is proved.

Moreover due to \eqref{MONo} and \eqref{limpun},
we have that for every fixed $\sigma>0$,
the sets
$${\mathcal A}_k^\sigma \equiv \{x\in {\mathbf R}^n| f_k(x)>\sigma\},$$
satisfy the following properties:
\begin{equation}\label{zapvar}
{\mathcal A}_{k+1}^\sigma \subset 
{\mathcal A}_{k}^\sigma \hbox{ } \forall k\in {\mathbf N}\hbox{ } 
\forall \sigma \in {\mathbf R}^+ 
\hbox{ and } meas (\cap_{k\in \mathbf N} {\mathcal A}_k^\sigma)=0
\hbox{  } \forall \sigma \in {\mathbf R}^+.
\end{equation}

On the other hand, since $f_1\in L^{p, d}(\Omega^n)$
with $p, d\neq \infty$, it is easy to deduce that
$meas ({\mathcal A}_1^\sigma)<\infty$ for every $\sigma>0$. By combining this fact with \eqref{zapvar}, 
we can deduce that:
$$\lim_{k\rightarrow \infty} m(\sigma, f_k)=\lim_{k\rightarrow \infty} meas({\mathcal A}_k^\sigma)
=0 
\hbox{  } \forall \sigma>0.$$

In particular for every $\epsilon>0$ there exists $k(\epsilon)\in \mathbf N$ such
that $$m(\epsilon, f_{k(\epsilon)})<\epsilon.$$

This inequality implies that
if $t>0$ is a fixed number, then
$$f_{k(\epsilon)}^*(t)\equiv inf\{\sigma | m(\sigma, f_{k(\epsilon)})<t\} < \epsilon,$$
provided that $0<\epsilon <t$.

This estimate, combined with the monotonicity of
$\{f_k^*(t)\}_{k\in \mathbf N}$ (see \eqref{cuppe}), 
implies easily \eqref{null}.

\hfill$\Box$

\noindent {\bf Proof of proposition \ref{propr}.}
Let us introduce the following
sets:
$$\Omega_k\equiv \{x\in \Omega | |a(x)|<k \hbox{ and } |x|<k\},$$
where $k\in \mathbf N$.

It is easy to show that the sequence of sets 
$\{\Omega_k\}_{k\in \mathbf N}$ satisfy the
following conditions:
$$\Omega_k \subset \Omega_{k+1}, |\Omega_k|<\infty \hbox{  } \forall k\in {\mathbf N}
\hbox{ and } a(x)\chi_{\Omega_k}\in L^\infty(\Omega).$$

It is then sufficient to prove that for every fixed $\epsilon>0$,
there exists of a suitable $k_0(\epsilon)\in \mathbf N$
such that 
\begin{equation}\label{final}
\|a(x) \chi_{\Omega \setminus 
\Omega_{k_0(\epsilon)}}\|_{L^{p, d}(\Omega)}<\epsilon,
\end{equation}
in order to conclude that the sets
$$\Omega^{\epsilon}_1=\Omega_{k_0(\epsilon)} \hbox { and }
\Omega^{\epsilon}_2=\Omega\setminus \Omega_{k_0(\epsilon)},$$
satisfy the desired properties.
 
In order to show \eqref{final} let us introduce the sequence
of functions $$a_k^*(t):{\mathbf R}^+ \rightarrow {\mathbf R}^+$$
where:
$$a_k^*(t)=(|a(x)|\chi_{\Omega \setminus \Omega_k})^*(t) \hbox{  } \forall k\in \mathbf N,$$
and as usual $f^*$ denotes
the decreasing rearrangement of the  
function $f$.

If $|a|^*(t)$ denotes
the rearranged function associated to $|a|$, 
then  by using lemma \ref{monotono} we deduce that
the sequence $\{a_k^*(t)\}_{k\in \mathbf N}$
satisfies the following inequalities:

\begin{align}
\label{1}&|a|^*(t)\geq a_k^*(t)\geq 0 \hbox{  } 
\forall t>0 \hbox{ and } \forall k\in {\mathbf N},\\
\label{3}&\lim_{k\rightarrow \infty} a_k^*(t)=0 \hbox{  } \forall t\in 
{\mathbf R}^+.
\end{align}

In particular, since $a(x)\in L^{p, d}(\Omega)$, we have
$$|a_k^*(t)|^d  t^{\frac{d}{p}-1}\leq 
(|a|^*(t))^d  t^{\frac{d}{p}-1} \in L^1({\mathbf R}^+) 
$$ 
that can be combined with the dominated convergence theorem
and with \eqref{3} in order to deduce that
$$\lim_{k\rightarrow \infty}
\|a(x)\chi_{\Omega \setminus \Omega_k} 
\|_{L^{p, d}(\Omega)}^d=\lim_{k\rightarrow \infty} \int_0^\infty |a_k^*(t)|^d  t^{\frac{d}{p}-1}
dt=0.$$

The proof is complete.

\hfill$\Box$

\end{document}